\documentclass{article}

\usepackage{amsmath,amssymb,amsthm}

\newtheorem{thm}{Theorem}
\newtheorem{defn}[thm]{Definition}
\newtheorem{rem}[thm]{Remark}

\begin{document}

\title{Conservation laws for invariant functionals
containing compositions\footnote{Accepted for an oral presentation
at the 7th IFAC Symposium on Nonlinear Control Systems (NOLCOS
2007), to be held in Pretoria, South Africa, 22--24 August,
2007.}}

\author{Gast\~{a}o S. F. Frederico\thanks{This work is part
         of the author's PhD project. Supported by the
         \emph{Portuguese Institute for Development} (IPAD).}\\
        \texttt{gfrederico@mat.ua.pt}\\[0.3cm]
        Higher Institute of Education\\
        University of Cabo Verde\\
        Praia, Santiago -- Cape Verde \and
        Delfim F. M. Torres\thanks{Supported by the
          \emph{Centre for Research on Optimization and Control} (CEOC)
          through the \emph{Portuguese Foundation for Science and Technology}
          (FCT), cofinanced by the European Community fund FEDER/POCTI.}\\
        \texttt{delfim@mat.ua.pt}\\[0.3cm]
        Department of Mathematics\\
        University of Aveiro\\
        3810-193 Aveiro, Portugal}

\date{}

\maketitle


\begin{abstract}
The study of problems of the calculus of variations with
compositions is a quite recent subject with origin in dynamical
systems governed by chaotic maps. Available results are reduced to
a generalized Euler-Lagrange equation that contains a new term
involving inverse images of the minimizing trajectories. In this
work we prove a generalization of the necessary optimality
condition of DuBois-Reymond for variational problems with
compositions. With the help of the new obtained condition, a
Noether-type theorem is proved. An application of our main result
is given to a problem appearing in the chaotic setting when one
consider maps that are ergodic.
\end{abstract}

\smallskip

\noindent \textbf{Mathematics Subject Classification 2000:} 49K05, 49J05.

\smallskip


\smallskip

\noindent \textbf{Keywords:} variational calculus, functionals containing compositions,\\
symmetries, DuBois-Reymond condition, Noether's theorem.

\medskip


\section{Introduction and motivation}

The theory of variational calculus for problems with compositions
has been recently initiated in \cite{CD:Bracken:2004}. The new
theory considers integral functionals that depend not only on
functions $q(\cdot)$ and their derivatives $\dot{q}(\cdot)$, but
also on compositions $(q \circ q)(\cdot)$ of $q(\cdot)$ with
$q(\cdot)$. As far as chaos is often a byproduct of iteration of
nonlinear maps \cite{CD:BraGo:1997}, such problems serve as an
interesting model for chaotic dynamical systems. Let us briefly
review this relation (for more details, we refer the interested
reader to \cite{CD:BraGoBo1:2001,CD:BraGoBo2:2002,CD:Bracken:2004}).
Let $q : [0,1] \rightarrow [0,1]$ be a piecewise monotonic map with
probability density function $f_q(\cdot)$, which captures the long
term statistical behavior of a nonlinear dynamical system. It is
natural (see \cite{CD:BraGo:1997}) to consider the
problem of minimizing or maximizing the functional
\begin{equation}
\label{eq:prbChaos}
I[q(\cdot),f_q(\cdot)] = \int_0^1 \left(q(t) - t\right)^2 f_q(t) dt \, ,
\end{equation}
which depends on $q(\cdot)$ and its probability density function
$f_q(\cdot)$ (usually a complicated function of $q(\cdot)$).
It turns out that $f_q(\cdot)$ is the fixed point of
the Frobenius-Perron operator $P_q[\cdot]$ associated with
$q(\cdot)$. For a piecewise monotonic map $q : [0,1] \rightarrow
[0,1]$ with $r$ pieces, $P_q[\cdot]$ has the representation
\begin{equation*}
P_q[f](t) = \sum_{v \in \{q^{-1}(t)\}}
\frac{f(v)}{\left|\dot{q}(v)\right|} \, ,
\end{equation*}
where for any point $t \in [0,1]$ the set $\{q^{-1}(t)\}$ consists
of at most $r$ points. The fixed point $f_q(\cdot)$ associated
with an ergodic map $q(\cdot)$ can be expressed as the limit
\begin{equation}
\label{eq:fq}
f_q = \lim_{n \rightarrow \infty} \sum_{i=0}^{n-1}
P_q^i[\mathbf{1}] \, ,
\end{equation}
where $\mathbf{1}$ is the constant function $1$ on $[0,1]$.
Substituting \eqref{eq:fq} into \eqref{eq:prbChaos}, and using the
adjoint property \cite[Prop.~4.2.6]{CD:BraGo:1997},
one eliminates the probability density function $f_q(\cdot)$,
obtaining \eqref{eq:prbChaos} in the form
\begin{equation*}
I[q(\cdot)] = \int_0^1
L\left(t,q(t),q^{(2)}(t),q^{(3)}(t),\ldots\right) dt \, ,
\end{equation*}
where we are using the notation $q^{(i)}(\cdot)$ to denote the
$i$-th composition of $q(\cdot)$ with itself: $q^{(1)}(t) = q(t)$,
$q^{(2)}(t) = (q \circ q)(t)$, $q^{(3)}(t) = (q \circ q \circ
q)(t)$, etc. In \cite{CD:Bracken:2004} a generalized
Euler-Lagrange equation, which involves the inverse images of the
extremizing function $q(\cdot)$ (\textrm{cf.} \eqref{eq:elfc}),
was proved for such functionals in the cases
$$\int_a^b L\left(t,q(t),q^{(2)}(t)\right) dt \, ,$$
$$\int_a^b L\left(t,q(t),\dot{q}(t),q^{(2)}(t)\right) dt\, ,$$
or
$$\int_a^b L\left(t,q(t),q^{(2)}(t),q^{(3)}(t)\right) dt\, .$$
To the best of our knowledge, these generalized Euler-Lagrange
equations comprise all the available results on the subject.
Thus, one concludes that the theory of variational calculus
with compositions is in its childhood: much remains to be done.
Here we go a step further in
the theory of functionals containing compositions. We are mainly
interested in Noether's classical theorem, which is one of the
most beautiful results of the calculus of variations and optimal control,
with many important applications in Physics (see \textrm{e.g.}
\cite{CD:Djukic:1980,CD:Logan:1987,MR1694555}), Economics
(see \textrm{e.g.} \cite{Askenazy,Sato}),
and Control Engineering (see \textrm{e.g.}
\cite{Gugushvili,MR83k:93011,MR96i:49037,delfim3ncnw,MR83k:49054}),
and source of many recent extensions and developments (see
\textrm{e.g.} \cite{CD:Gastao:2006,GastaoIJTS,CD:LiQun:2003,%
GouveiaTorresCMAM,comEugenioCC,delfimEJC,CD:JMS:Torres:2004}).
Noether's symmetry theorem describes the universal fact
that invariance with respect to some family of parameter
transformations gives rise to the existence of certain conservation laws,
\textrm{i.e.} expressions preserved along the Euler-Lagrange
or Pontryagin extremals of the problem.
Our results are a generalized DuBois-Reymond necessary optimality
condition (Theorem~\ref{theo:cDRfc}), and a generalized Noether's
theorem (Theorem~\ref{theo:tnfc}) for functionals of the form
$\int_a^b L\left(t,q(t),\dot{q}(t),q^{(2)}(t)\right) dt$. In
\S\ref{sec:ex} an illustrative example is presented.


\section{Preliminaries -- review of classical results of the calculus of variations}

There exist many different ways to prove the classical Noether's
theorem (\textrm{cf. e.g.}
\cite{CD:Djukic:1980,CD:Jost:1998,CD:Logan:1987,Sato}). We review
here one of those proofs, which is based on the DuBois-Reymond
necessary condition. Although this proof is not so common in the
literature of Noether's theorem, it turns out to be the most
suitable approach when dealing with functionals containing
compositions.

Let us consider the fundamental problem of the calculus of
variations:
\begin{equation}
\label{P}
I[q(\cdot)] = \int_a^b L\left(t,q(t),\dot{q}(t)\right) dt
\longrightarrow \min  \tag{P}
\end{equation}
under the boundary conditions $q(a)=q_{a}$ and $q(b)=q_{b}$, where
$\dot{q} = \frac{dq}{dt}$, with $q(\cdot)$ a piecewise-smooth function,
and the Lagrangian $L :[a,b] \times
\mathbb{R}^{n} \times \mathbb{R}^{n} \rightarrow \mathbb{R}$ is a
$C^{2}$ function with respect to all its arguments.

The concept of symmetry has a very important role in mathematics
and its applications. Symmetries are defined through
transformations of the system that leave the problem
\emph{invariant}.

\begin{defn}[Invariance of \eqref{P}]
\label{def:inva} The integral functional \eqref{P} it said to be
invariant under the $\varepsilon$-parameter
infinitesimal transformations
\begin{equation}
\label{eq:tinf}
\begin{cases}
\bar{t} = t + \varepsilon\tau(t,q) + o(\varepsilon) \, ,\\
\bar{q}(t) = q(t) + \varepsilon\xi(t,q) + o(\varepsilon) \, ,\\
\end{cases}
\end{equation}
where $\tau$ and $\xi$ are piecewise-smooth, if
\begin{equation}
\label{eq:inv1} \int_{t_{a}}^{t_{b}}
L\left(t,q(t),\dot{q}(t)\right)dt =
\int_{\bar{t}(t_a)}^{\bar{t}(t_b)}
L\left(\bar{t},\bar{q}(\bar{t}),\dot{\bar{q}}(\bar{t})\right)d\bar{t}
\end{equation}
for any subinterval $[{t_{a}},{t_{b}}] \subseteq [a,b]$.
\end{defn}

Along the work we denote by $\partial_{i}L$ the partial derivative
of $L$ with respect to its $i$-th argument.

\begin{thm}[Necessary condition of invariance]
\label{theo:cnsi} If functional \eqref{P} is invariant under the
infinitesimal transformations \eqref{eq:tinf}, then
\begin{equation}
\label{eq:cnsi}
\partial_{1}
L\left(t,q,\dot{q}\right)\tau+\partial_{2}
L\left(t,q,\dot{q}\right)\cdot\xi
+ \partial_{3} L\left(t,q,\dot{q}\right)\cdot
\left(\dot{\xi}-\dot{q}\dot{\tau}\right)
+ L\left(t,q,\dot{q}\right)\dot{\tau}=0 \, .
\end{equation}
\end{thm}

\begin{proof}
Since \eqref{eq:inv1} is to be satisfied for any subinterval
$[{t_{a}},{t_{b}}] \subseteq [a,b]$, equality \eqref{eq:inv1} is
equivalent to
\begin{equation}
\label{eq:inv2}
\left[ L\left(t+\varepsilon\tau+o(\varepsilon),q+\varepsilon\xi
+o(\varepsilon), \frac{\dot{q}+ \varepsilon\dot{\xi}
+o(\varepsilon)}{1+\varepsilon \dot{\tau}
+o(\varepsilon)}\right)\right]\frac{d\bar{t}}{dt}
= L\left(t,q,\dot{q}\right)\, .
\end{equation}
We obtain \eqref{eq:cnsi} differentiating both sides of
\eqref{eq:inv2} with respect to $\varepsilon$, and then setting
$\varepsilon=0$.
\end{proof}

Another very important notion in mathematics and its applications
is the concept of \emph{conservation law}. One of the most
important conservation laws was proved by Leonhard Euler in 1744:
when the Lagrangian $L(q,\dot{q})$ corresponds to a system of
conservative points, then
\begin{equation}
\label{eq:consEneg} -L\left(q(t),\dot{q}(t)\right) +\frac{\partial
L}{\partial \dot{q}}\left(q(t),\dot{q}(t)\right) \cdot \dot{q}(t)
\equiv \text{constant} \, ,
\end{equation}
$t \in [a,b]$, holds along the solutions
of the Euler-Lagrange equations.

\begin{defn}[Conservation law]
\label{def:leico} A quantity $C(t,q,\dot{q})$
defines a conservation law if
$$\frac{d}{dt}C(t,q(t),\dot{q}(t))=0 \, ,
\quad t \in [a,b] \, ,$$
along all the solutions $q(\cdot)$ of the Euler-Lagrange equation
\begin{equation}
\label{eq:el} \frac{d}{{dt}}\partial_{3} L\left(t,q,\dot{q}\right) =
\partial_{2} L\left(t,q,\dot{q}\right) \, .
\end{equation}
\end{defn}

Conservation laws can be used to lower the order of the
Euler-Lagrange equations \eqref{eq:el} and simplify the
resolution of the respective problems of the calculus of
variations and optimal control \cite{comEugenioCC}. Emmy Amalie
Noether formulated in 1918 a very general principle on
conservation laws, with many important implications in modern
physics, economics and engineering. Noether's principle asserts that
\emph{``the invariance of the functional $\int_a^b
L\left(t,q(t),\dot{q}(t)\right) dt$ under one-parameter
infinitesimal transformations \eqref{eq:tinf}, imply the existence
of a conservation law''}. One particular example of application of
Noether's theorem gives \eqref{eq:consEneg}, which corresponds to
conservation of energy in classical mechanics or to the
income-wealth law of economics.

\begin{thm}[Noether's theorem]
\label{theo:tnoe} If functional \eqref{P} is invariant, in the
sense of the Definition~\ref{def:inva}, then
\begin{equation}
\label{eq:TeNet} C(t,q,\dot{q}) =
\partial_{3} L\left(t,q,\dot{q}\right)\cdot\xi(t,q)
+ \left( L(t,q,\dot{q}) - \partial_{3} L\left(t,q,\dot{q}\right)
\cdot \dot{q} \right) \tau(t,q)
\end{equation}
defines a conservation law.
\end{thm}

We recall here the proof of Theorem~\ref{theo:tnoe} by means of
the classical necessary optimality condition of DuBois-Reymond.

\begin{thm}[DuBois-Reymond condition]
\label{theo:cdr} If $q(\cdot)$ is a solution of problem \eqref{P},
then
\begin{equation}
\label{eq:cdr}
\partial_{1}
L\left(t,q,\dot{q}\right)=\frac{d}{dt}\left[L\left(t,q,\dot{q}\right)-\partial_{3}
L\left(t,q,\dot{q}\right)\cdot\dot{q}\right]\, .
\end{equation}
\end{thm}

\begin{proof}
The DuBois-Reymond necessary optimality
condition is easily proved using the
Euler-Lagrange equation \eqref{eq:el}:
\begin{equation*}
\begin{split}
\frac{d}{dt}&\left[L\left(t,q,\dot{q}\right)-\partial_{3}
L\left(t,q,\dot{q}\right)\cdot\dot{q}\right]\\
& =\partial_{1}
L\left(t,q,\dot{q}\right)+\partial_{2}
L\left(t,q,\dot{q}\right)\cdot\dot{q} +\partial_{3}
L\left(t,q,\dot{q}\right)\cdot\ddot{q} \\
& \quad -\frac{d}{dt}\partial_{3}
L\left(t,q,\dot{q}\right)\cdot\dot{q}-\partial_{3}
L\left(t,q,\dot{q}\right)\cdot\ddot{q} \\
& =\partial_{1}
L\left(t,q,\dot{q}\right)+\dot{q}\cdot\left(\partial_{2}
L\left(t,q,\dot{q}\right)-\frac{d}{dt}\partial_{3}
L\left(t,q,\dot{q}\right)\right)\\
&=\partial_{1} L\left(t,q,\dot{q}\right)\, .
\end{split}
\end{equation*}
\end{proof}

\begin{proof} (of Theorem~\ref{theo:tnoe})
To prove the Noether's theorem, we use the
Euler-Lagrange equation \eqref{eq:el} and the DuBois-Reymond condition
\eqref{eq:cdr} into the necessary condition of invariance
\eqref{eq:cnsi}:
\begin{equation*}
\begin{split}
0 &= \partial_{1} L\left(t,q,\dot{q}\right)\tau
 +\partial_{2} L\left(t,q,\dot{q}\right)\cdot\xi \\
& \quad +\partial_{3} L\left(t,q,\dot{q}\right)\cdot
\left(\dot{\xi}-\dot{q}\dot{\tau}\right)
+L\left(t,q,\dot{q}\right)\dot{\tau} \\
&=\partial_{2} L\left(t,q,\dot{q}\right)\cdot\xi+\partial_{3}
L\left(t,q,\dot{q}\right)\cdot\dot{\xi} +\partial_{1}
L\left(t,q,\dot{q}\right)\tau\\
& \quad +\dot{\tau}\left(L\left(t,q,\dot{q}\right)-\partial_{3}
L\left(t,q,\dot{q}\right)\cdot\dot{q}\right)\\
&=\frac{d}{{dt}}\partial_{3}
L\left(t,q,\dot{q}\right)\cdot\xi+\partial_{3}
L\left(t,q,\dot{q}\right)\cdot\dot{\xi}\\
& \quad +\frac{d}{dt}\left(L\left(t,q,\dot{q}\right)
-\partial_{3} L\left(t,q,\dot{q}\right)\cdot\dot{q}\right)\tau \\
& \quad +\dot{\tau}\left(L\left(t,q,\dot{q}\right)-\partial_{3}
L\left(t,q,\dot{q}\right)\cdot\dot{q}\right)\\
&=\frac{d}{dt}\Bigl[\partial_{3} L\left(t,q,\dot{q}\right)\cdot\xi
+ \bigl( L(t,q,\dot{q}) - \partial_{3} L\left(t,q,\dot{q}\right) \cdot \dot{q} \bigr)
\tau\Bigr]\, .
\end{split}
\end{equation*}
\end{proof}


\section{Main results}
\label{sec:MR}

We consider the following problem of the calculus of
variations with composition of functions:
\begin{equation}
\label{Pc} I[q(\cdot)] = \int_a^b
L\left(t,q(t),\dot{q}(t),z(t)\right)dt \longrightarrow \min
\tag{$P_c$}
\end{equation}
subject to given boundary conditions $q(a)=q_{a}$, $q(b)=q_{b}$,
$z(a)=z_{a}$, and $z(b)=z_{b}$, where
$\dot{q} = \frac{dq}{dt}$ and $z(t)=(q\circ q)(t)$. We assume that
the Lagrangian $L :[a,b] \times \mathbb{R} \times \mathbb{R} \times \mathbb{R}
\rightarrow \mathbb{R}$ is a function of class $C^{2}$
with respect to all the arguments, and that admissible functions
$q(\cdot)$ are piecewise-smooth. The main result of \cite{CD:Bracken:2004}
is an extension of the Euler-Lagrange equation \eqref{eq:el}
for problems of the calculus of variations \eqref{Pc}.

\begin{thm}[\cite{CD:Bracken:2004}]
\label{Thm:E-Lfc} If $q(\cdot)$ is a weak minimizer of problem \eqref{Pc},
then $q(\cdot)$ satisfies the \emph{Euler-Lagrange equation}
\begin{multline}
\label{eq:elfc}
\partial_{2} L\left(x,q(x),\dot{q}(x),z(x)\right)
-\frac{d}{dx} \partial_{3} L\left(x,q(x),\dot{q}(x),z(x)\right)\\
+\partial_{4} L\left(x,q(x),\dot{q}(x),z(x)\right) \dot{q}(q(x))
+\sum_{t=q^{-1}(x)}\frac{\partial_{4}
L\left(t,q(t),\dot{q}(t),z(t)\right)}{|\dot{q}(t)|}=0
\end{multline}
for any $x \in (a,b)$.
\end{thm}


\subsection{Generalized DuBois-Reymond condition}

We begin by proving an extension of the
DuBois-Reymond necessary optimality condition \eqref{eq:cdr}
for problems of the calculus of variations \eqref{Pc}.

\begin{thm}[\textrm{cf.} Theorem~\ref{theo:cdr}]
\label{theo:cDRfc}
If $q(\cdot)$ is a weak minimizer of problem \eqref{Pc},
then $q(\cdot)$ satisfies the \emph{DuBois-Reymond condition}
\begin{multline}
\label{eq:cDRfc}
\frac{d}{dx}\Bigl[L\left(x,q(x),\dot{q}(x),z(x)\right)
-\partial_{3} L\left(x,q(x),\dot{q}(x),z(x)\right)\dot{q}(x)\Bigr]\\
=\partial_{1}
L\left(x,q(x),\dot{q}(x),z(x)\right)
-\dot{q}(x)\sum_{t=q^{-1}(x)}\frac{\partial_{4}
L\left(t,q(t),\dot{q}(t),z(t)\right)}{|\dot{q}(t)|}
\end{multline}
for any $x \in (a,b)$.
\end{thm}

\begin{rem}
If $L\left(t,q,\dot{q},z\right) = L\left(t,q,\dot{q}\right)$,
then \eqref{eq:cDRfc} coincides with the classical
DuBois-Reymond condition \eqref{eq:cdr}.
\end{rem}

\begin{proof}
To prove Theorem~\ref{theo:cDRfc} we use the
Euler-Lagrange equation \eqref{eq:elfc}:
\begin{equation*}
\begin{split}
\frac{d}{dx} & \Bigl[L\left(x,q,\dot{q},z\right)
- \partial_{3} L\left(x,q,\dot{q},z\right)\dot{q}\Bigr]\\
&=\partial_{1} L\left(x,q,\dot{q},z\right)
+\partial_{2} L\left(x,q,\dot{q},z\right)\dot{q}\\
& \quad +\partial_{3} L\left(x,q,\dot{q},z\right)\ddot{q}
+\partial_{4} L\left(x,q,\dot{q},z\right)\dot{q}(q(x))\dot{q}\\
& \quad -\dot{q}\frac{d}{dx}\partial_{3} L\left(x,q,\dot{q},z\right)
-\partial_{3} L\left(x,q,\dot{q},z\right)\ddot{q} \\
&=\partial_{1} L\left(x,q,\dot{q},z\right)
+\dot{q}\Bigl(\partial_{2} L\left(x,q,\dot{q},z\right)\\
& \quad + \partial_{4} L\left(x,q,\dot{q},z\right)\dot{q}(q(x))-\frac{d}{dx}\partial_{3}
L\left(x,q,\dot{q},z\right)\Bigr)\\
&=\partial_{1} L\left(x,q,\dot{q},z\right)
-\dot{q}(x)\sum_{t=q^{-1}(x)}\frac{\partial_{4}
L\left(t,q(t),\dot{q}(t),z(t)\right)}{|\dot{q}(t)|}\, .
\end{split}
\end{equation*}
\end{proof}


\subsection{Noether's theorem for functionals containing compositions}

We introduce now the definition of \emph{invariance}
for the functional \eqref{Pc}. As done
in the proof of Theorem~\ref{theo:cnsi} (see \eqref{eq:inv2}),
we get rid off of the integral signs in \eqref{eq:inv1}.

\begin{defn}[\textrm{cf.} Definition~\ref{def:inva}]
\label{def:invafc} We say that functional \eqref{Pc} is invariant
under the infinitesimal transformations \eqref{eq:tinf} if
\begin{equation}
\label{eq:invfc}
L\left(\bar{t},\bar{q}(\bar{t}),{\bar{q}}'(\bar{t}),\bar{z}(\bar{t})\right)
\frac{d\bar{t}}{d t}
= L\left(t,q(t),\dot{q}(t),z(t)\right)+o(\varepsilon) \, ,
\end{equation}
where $\bar{q}'={d\bar{q}}/{d\bar{t}}$.
\end{defn}

Along the work, in order to simplify the presentation,
we sometimes omit the arguments of the functions.

\begin{thm}[\textrm{cf.} Theorem~\ref{theo:cnsi}]
\label{theo:cnsifc} If functional \eqref{Pc} is invariant
under the infinitesimal transformations \eqref{eq:tinf}, then
\begin{multline}
\label{eq:cnsifc}
\partial_{1} L\left(t,q,\dot{q},z\right)\tau
+\partial_{2} L\left(t,q,\dot{q},z\right)\xi
+\partial_{3} L\left(t,q,\dot{q},z\right)\left(\dot{\xi}-\dot{q}\dot{\tau}\right)\\
+\partial_{4} L\left(t,q,\dot{q},z\right)\dot{q}(q(t))\xi
+ \partial_{4} L\left(t,q,\dot{q},z\right)\xi(q(t))+L\dot{\tau} =0\, .
\end{multline}
\end{thm}

\begin{proof}
Equation \eqref{eq:invfc} is equivalent to
\begin{multline}
\label{eq:invif2}
L\Biggl(t+\varepsilon\tau+o(\varepsilon),q+\varepsilon\xi+o(\varepsilon),
\frac{\dot{q}+ \varepsilon\dot{\xi}+o(\varepsilon)}{1+\varepsilon
\dot{\tau}+o(\varepsilon)},\\
q(q+\varepsilon\xi+o(\varepsilon))
+\varepsilon\xi(q+\varepsilon\xi+o(\varepsilon))\Biggr) \times
\left(1+\varepsilon \dot{\tau}+o(\varepsilon)\right) \\
= L\left(t,q,\dot{q},z\right) +o(\varepsilon) \, .
\end{multline}
We obtain equation \eqref{eq:cnsifc}
differentiating both sides of equality \eqref{eq:invif2}
with respect to the parameter $\varepsilon$, and then setting
$\varepsilon=0$.
\end{proof}

\begin{rem}
Using the Frobenius-Perron operator (see \cite[Chap.~4]{CD:BraGo:1997})
and the Euler-Lagrange equation \eqref{eq:elfc}, we can write
\eqref{eq:cnsifc} in the following form:
\begin{multline}
\label{eq:cnsifc1}
\partial_{1}
L\left(x,q,\dot{q},z\right)\tau+\partial_{2}
L\left(x,q,\dot{q},z\right)\xi \\
+\partial_{3} L\left(x,q,\dot{q},z\right)\left(\dot{\xi}
-\dot{q}\dot{\tau}\right)
+\partial_{4} L\left(x,q,\dot{q},z\right)\dot{q}(q(x))\xi\\
+ \sum_{t=q^{-1}(x)}\frac{\partial_{4}
L\left(t,q(t),\dot{q}(t),z(t)\right)}{|\dot{q}(t)|}\xi+L\dot{\tau}\\
=\partial_{1} L\left(x,q,\dot{q},z\right)\tau
+\frac{d}{dx}\partial_{3} L\left(x,q,\dot{q},z\right)\xi\\
+\partial_{3} L\left(x,q,\dot{q},z\right)\left(\dot{\xi}
-\dot{q}\dot{\tau}\right)+L\dot{\tau}=0\, .
\end{multline}
\end{rem}

\begin{defn}[Conservation law for \eqref{Pc}]
\label{def:leicofc} We say that a quantity
$C\left(x,q,\dot{q},z\right)$ defines
a \emph{conservation law for functionals
containing compositions} if
$$\frac{d}{dx}C\left(x,q(x),\dot{q}(x),z(x)\right)=0$$
along all the solutions $q(\cdot)$ of the
Euler-Lagrange equation \eqref{eq:elfc}.
\end{defn}

Our main result is an extension of Noether's theorem for
problems of the calculus of variations \eqref{Pc}
containing compositions.

\begin{thm}[Noether's theorem for \eqref{Pc}]
\label{theo:tnfc} If functional \eqref{Pc} is invariant, in the sense
of the Definition~\ref{def:invafc}, and there exists a function
$f=f(x,q,\dot{q},z)$ such that
\begin{equation}
\label{eq:condfc}
\frac{df}{dx}\left(x,q(x),\dot{q}(x),z(x)\right)
=\tau\dot{q}(x)\sum_{t=q^{-1}(x)}\frac{\partial_{4}
L\left(t,q(t),\dot{q}(t),z(t)\right)}{|\dot{q}(t)|} \, ,
\end{equation}
then
\begin{multline}
\label{eq:TeNetfc}
C\left(x,q(x),\dot{q}(x),z(x)\right) \\
= \Bigl[ L(x,q(x),\dot{q}(x),z(x))
- \partial_{3} L\left(x,q(x),\dot{q}(x),z(x)\right)
\dot{q} \Bigr] \tau(x,q)\\
+\partial_{3} L\left(x,q(x),\dot{q}(x),z(x)\right)\xi(x,q)
+f(x,q(x),\dot{q}(x),z(x))
\end{multline}
defines a conservation law (Definition~\ref{def:leicofc}).
\end{thm}

\begin{rem}
If $L\left(x,q,\dot{q},z\right)=L\left(x,q,\dot{q}\right)$, then
$f$ is a constant and expression \eqref{eq:TeNetfc} is equivalent
to the conserved quantity \eqref{eq:TeNet}
given by the classical Noether's theorem.
\end{rem}

\begin{proof}
To prove the theorem, we use conditions
\eqref{eq:cDRfc} and \eqref{eq:condfc} in \eqref{eq:cnsifc1}:
\begin{equation*}
\begin{split}
0 &= \partial_{1} L\left(x,q,\dot{q},z\right)\tau
+ \frac{d}{dx}\partial_{3} L\left(x,q,\dot{q},z\right)\xi\\
& \quad +\partial_{3} L\left(x,q,\dot{q},z\right)
\left(\dot{\xi}-\dot{q}\dot{\tau}\right)+L\dot{\tau}\\
&=\tau\frac{d}{dx}\Bigl[L\left(x,q,\dot{q},z\right)-\partial_{3}
L\left(x,q,\dot{q}),z)\right)\dot{q}\Bigr]\\
& \quad + \dot{\tau}\left[L\left(x,q,\dot{q},z\right)-\partial_{3}
L\left(x,q,\dot{q}),z)\right)\dot{q}\right]\\
& \quad +\dot{\xi}\partial_{3}
L\left(x,q,\dot{q},z\right)+\xi\frac{d}{dx}\partial_{3}
L\left(x,q,\dot{q},z\right)\\
& \quad +\tau \dot{q}(x)\sum_{t=q^{-1}(x)}\frac{\partial_{4}
L\left(t,q(t),\dot{q}(t),z(t)\right)}{|\dot{q}(t)|}\\
&=\frac{d}{dx}\Biggl\{\partial_{3}
L\left(x,q,\dot{q},z\right)\xi
+\Bigl[L\left(x,q,\dot{q},z\right)\\
& \quad -\partial_{3}
L\left(x,q,\dot{q},z\right)\dot{q}\Bigr]\tau
+f(x,q,\dot{q},z)\Biggr\} \, .
\end{split}
\end{equation*}
\end{proof}


\section{An example}
\label{sec:ex}

Let us consider the problem
\begin{gather}
I[q(\cdot)] = \frac{1}{3}\int_0^1 \left[x+q(x)+q(q(x))\right] dx
\longrightarrow \min \notag \\
q(0) = 1 \, , \quad q(1) = 0 \, , \label{eq:exemp} \\
q(q(0)) = 0 \, , \quad q(q(1)) = 1 \, . \notag
\end{gather}
In \cite[\S 3]{CD:Bracken:2004} it is proven
that \eqref{eq:exemp} has the extremal
\begin{equation}
\label{eq:extre1}
q(x)=\begin{cases}
q_1(x)=-2x+1\, , \qquad x\in\left[0,\frac{1}{2}\right[ \, , \\
q_2(x)=-2x+2\, , \qquad x\in\left[\frac{1}{2},1\right]\, ,
\end{cases}
\end{equation}
that is, \eqref{eq:extre1} satisfies the Euler-Lagrange equation
\eqref{eq:elfc} for $L(x,q,\dot{q},z)=\frac{1}{3}\left(x+q+z\right)$.
We now illustrate the application of our Theorem~\ref{theo:tnfc}
to this problem. First, we need to determine the variational symmetries.
Substituting the Lagrangian $L$ in \eqref{eq:cnsifc1}
we obtain that
\begin{equation}
\label{eq:tau}
\frac{\tau}{3}+\frac{x+q+z}{3}\dot{\tau}=0\, .
\end{equation}
The differential equation \eqref{eq:tau} admits the solution
\begin{equation}
\label{eq:tau1}
 \tau=ke^{-\int\frac{dx}{x+q+z}} \, ,
\end{equation}
where $k$ is an arbitrary constant. From Theorem~\ref{theo:tnfc}
we conclude that
\begin{equation}
\label{eq:con1}
 (x+q_{1}+z_1)\tau +\frac{1}{3}\int\tau
 \dot{q_1}\sum_{t=q^{-1}_1(x)}\frac{1}{|\dot{q_1}(t)|}dx, \quad
 x\in\left[0,\frac{1}{2}\right[ \, ,
\end{equation}
\begin{equation}
\label{eq:con2}
 (x+q_{2}+z_2)\tau
 +\frac{1}{3}\int\tau
 \dot{q_2}\sum_{t=q^{-1}_2(x)}\frac{1}{|\dot{q_2}(t)|}dx, \quad
 x\in\left[\frac{1}{2},1\right] \, ,
\end{equation}
defines a conservation law, where $\tau$ is obtained from \eqref{eq:tau1}:
\begin{equation}
\label{eq:tau2}
 \tau=ke^{-\int\frac{dx}{3x}}=ke^{\ln
 x^{-\frac{1}{3}}}=kx^{-\frac{1}{3}},
 \,\,  x\in[0,1]\, .
\end{equation}
Since for this problem we know the extremal, we can verify
the validity of the obtained conservation law
directly from Definition~\ref{def:leicofc}:
substituting equalities \eqref{eq:extre1} and \eqref{eq:tau2}
in \eqref{eq:con1} and \eqref{eq:con2},
we obtain, as expected, a constant (zero in this case):
\begin{multline*}
(x+q_{1}+z_1)\tau+\frac{1}{3}\int\tau
 \dot{q_1}\sum_{t=q^{-1}_1(x)}\frac{1}{|\dot{q_1}(t)|}dx\\
 =3kx\tau-2 \int\tau dx=3kx^{\frac{2}{3}}-3kx^{\frac{2}{3}}=0\, ,
\end{multline*}
\begin{multline*}
(x+q_{2}+z_2)\tau+\frac{1}{3}\int\tau
\dot{q_2}\sum_{t=q^{-1}_2(x)}\frac{1}{|\dot{q_2}(t)|}dx\\
= 3kx\tau-2 \int\tau dx=3kx^{\frac{2}{3}}-3kx^{\frac{2}{3}}=0\, .
\end{multline*}


\section{Conclusions}

We proved a generalization (i) of the necessary optimality
condition of DuBois-Reymond, (ii) of the celebrated
Noether's symmetry theorem, for problems
of the calculus of variations containing compositions
(respectively Theorems~\ref{theo:cDRfc} and \ref{theo:tnfc}).
Our main result is illustrated with the example studied
in \cite{CD:Bracken:2004}.

The compositional variational theory is in its childhood, so that
much remains to be done. In particular, it would be interesting to
obtain an Hamiltonian formulation and to study more general
optimal control problems with compositions.


\section*{Acknowledgements}

The authors are grateful to Pawe{\l} G{\'o}ra
who shared Chapter~4 of \cite{CD:BraGo:1997}.


\end{document}